\RequirePackage{amsmath}
\documentclass[a4paper,11pt]{article}

\usepackage{mathptmx}
\usepackage[english]{babel}

\usepackage{amsfonts}
\usepackage{amssymb}
\usepackage{amsthm}
\usepackage{graphicx}
\usepackage{pdfpages}
\usepackage{relsize}

\newcommand{\setx}[1]{ \{ #1 \} }

\newcommand{\mytext}[1]{ \: \textrm{#1} \: }
\newcommand{\mysetdescr}[2]{\left\{ #1 \: \left| \: #2 \right. \right\} }
\newcommand{\darr}{{\downarrow \,}}
\newcommand{\uarr}{{\uparrow \,}}

\newcommand\urbild[1]{#1^{-1}}

\def\yps{{\upsilon}}

\def\C{{\cal C}}

\def\F{{\cal F}}

\def\J{{\cal J}}
\def\K{{\cal K}}
\def\L{{\cal L}}
\def\M{{\cal M}}
\def\N{{\cal N}}

\def\U{{\cal U}}

\newcommand{\mf}[1]{\mathfrak{ #1 }}
\newcommand{\fC}{\mf{C}}
\newcommand{\fF}{\mf{F}}

\def\BP{\begin{proof}}
\def\EP{\end{proof}}

\DeclareMathOperator{\id}{id}
\DeclareMathOperator{\ImCr}{ImCr}
\DeclareMathOperator{\imp}{imp}

\newcommand{\CINull}{\ImCr}
\newcommand{\Co}{\CINull}
\newcommand{\Codr}{\CINull_{\max}}
\newcommand{\fCdr}{\fC_{\max}}

\newcommand{\FP}{\F(P)}
\newcommand{\fFP}{\fF(P)}
\newcommand{\fFQ}{\fF(Q)}
\newcommand{\MP}{\M(P)}

\newcommand{\Mvee}{{\displaystyle \mathlarger{\curlyvee}}}
\newcommand{\Mwedge}{{\displaystyle \mathlarger{\curlywedge}}}

\begin{document}

\theoremstyle{plain}
\theoremstyle{plain}
\newtheorem{theorem}{Theorem}
\newtheorem{definition}{Definition}
\newtheorem{corollary}{Corollary}
\newtheorem{lemma}{Lemma}
\newtheorem{proposition}{Proposition}
\newtheorem{result}{Result}

\title{{\bf About posets of height one as retracts}}
\author{\sc Frank a Campo}
\date{\small Viersen, Germany\\
{\sf acampo.frank@gmail.com}}

\maketitle
\begin{abstract}
\noindent We investigate connected posets $C$ of height one as retracts of finite posets $P$. We define two multigraphs: a multigraph $\fFP$ reflecting the network of so-called improper 4-crown bundles contained in the extremal points of $P$, and a multigraph $\fC(C)$ depending on $C$ but not on $P$. There exists a close interdependence between $C$ being a retract of $P$ and the existence of a graph homomorphism of a certain type from $\fFP$ to $\fC(C)$. In particular, if $C$ is an ordinal sum of two antichains, then $C$ is a retract of $P$ iff such a graph homomorphism exists. Returning to general connected posets $C$ of height one, we show that the image of such a graph homomorphism can be a clique in $\fC(C)$ iff the improper 4-crowns in $P$ contain only a sparse subset of the edges of $C$.
\newline

\noindent{\bf Mathematics Subject Classification:}\\
Primary: 06A07. Secondary: 06A06.\\[2mm]
{\bf Key words:} fixed point, fixed point property, retract, retraction, 4-crown, ordinal sum.
\end{abstract}

\section{Introduction}

A fundamental tool in the investigation of the fixed point property are {\em retractions} and {\em retracts}. These objects have been developed and intensively studied by a group of scientists around Rival \cite{Bjoerner_Rival_1980,Duffus_etal_1980_DPR,Duffus_Rival_1979,
Duffus_Rival_1981,Duffus_etal_1980_DRS,Nowakowski_Rival_1979,
Rival_1976,Rival_1980,Rival_1982} in the second half of the seventies. From the later work about the fixed point property and retracts, pars pro toto the work of Rutkowski \cite{Rutkowski_1986,Rutkowski_1989,Rutkowski_Schroeder_1994}, Schr\"{o}der \cite{Rutkowski_Schroeder_1994,Schroeder_1993,Schroeder_1995,Schroeder_1996,Schroeder_2012,Schroeder_2021,Schroeder_2022_MASoC}, and Zaguia \cite{Zaguia_2008} is mentioned here. For details, the reader is referred to \cite{Schroeder_2016,Schroeder_2012}.

In a previous paper \cite{aCampo_2024}, the author has characterized 4-crowns being retracts. Given a finite poset $P$, the approach based on a multigraph $\fFP$ reflecting the network of so-called 4-crown bundles contained in the extremal points $E(P)$ of $P$. The result was that a 4-crown in $E(P)$ is a retract of $P$ iff there exists a so-called $C$-separating homomorphism from $\fFP$ to a multigraph $\fC$ not depending on $P$.

This concept is generalized in Section \ref{subsec_Cgraph} of the present paper. Let $C$ be a connected poset of height one with carrier $Z$. We define a multigraph $\fC(C)$ whose vertices are the possible non-singleton images of 4-crown bundles under order-preserving mappings to $C$. The definition of a $Z$-separating homomorphism from $\fFP$ to $\fC(C)$ is adapted accordingly.


With these structures, we generalize the main results of \cite{aCampo_2024}. In Section \ref{subsec_sepHomsAndRetracts}, we see in Lemma \ref{lemma_surjAufC} and Lemma \ref{lemma_heightOne} that there is a close connection between $C \subseteq E(P)$ being a retract of $P$ and the existence of a $Z$-separating homomorphism from $\fFP$ to $\fC(C)$. In particular, we conclude in Theorem \ref{theo_surjAufC} that for $C$ being an ordinal sum of two antichains, $C$ is a retract of $P$ iff a $Z$-separating homomorphism $\phi : \fFP \rightarrow \fC(C)$ exists. Returning to general connected posets $C$ of height one, we prove in Section \ref{subsec_imagePhi} in Proposition \ref{prop_clique} that the image $\phi[\fFP]$ of a $Z$-separating homomorphism can be a clique in $\fC(C)$ iff the 4-crown bundles in $P$ contain only a sparse subset of the edges of $C$.



\section{Preparation} \label{sec_preparation}

\subsection{Notation and basic terms} \label{subsec_notation}

In what follows, $P = (X, \leq_P)$ is a finite poset. The pairs $x, y \in P$ with $x <_P y$ are called the {\em edges of $P$}. For $Y \subseteq X$, the {\em induced sub-poset} $P \vert_Y$ of $P$ is defined as $\left( Y, {\leq_P} \cap (Y \times Y) \right)$. To simplify notation, we identify a subset $Y \subseteq X$ with the poset $P \vert_Y$ induced by it.

For $y \in X$, we define the {\em down-set} and {\em up-set induced by $y$} as
\begin{align*}
\darr_P \; y & := \mysetdescr{ x \in X }{ x \leq_P y }, \quad
\uarr_P \; y := \mysetdescr{ x \in X }{ y \leq_P x },\end{align*}
and for $x, y \in X$, the {\em interval $[x,y]_P$} is defined as $[x,y]_P := ( \uarr_P x ) \cap ( \darr_P y )$.

We define the following subsets of $X$:
\begin{align*}
L(P) & := \mytext{the set of minimal points of } P, \\
U(P) & := \mytext{the set of maximal points of } P, \\
E(P) & := L(P) \cup U(P) \mytext{ the set of extremal points of } P, \\
M(P) & := X \setminus E(P).
\end{align*}
Following our convention to identify a subset of $X$ with the sub-poset induced by it, the minimal (maximal) points of $Y \subseteq X$ are the minimal (maximal) points of $P \vert_Y $, and we write $L(Y)$ instead of $L( P \vert_Y)$ and analogously for the other sets.

\begin{figure}
\begin{center}
\includegraphics[trim = 70 715 310 75, clip]{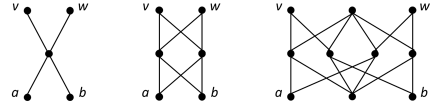}
\caption{\label{fig_Intro} Three posets $P$ of height two in which the points of $C = \setx{a,b,v,w}$ form a 4-crown in $E(P)$. In the first two posets, $C$ is an improper 4-crown. In the rightmost poset, $C$ is a proper 4-crown, and every edge of every 4-crown in $E(P)$ belongs to an improper 4-crown.}
\end{center}
\end{figure}

The following definition is crucial:

\begin{definition} \label{def_inner}
Let $C = \setx{a,b,v,w}$ be a 4-crown in $P$ with $\setx{a,b} = L(C)$ and $\setx{v,w} = U(C)$. We call $\J_P(C) := [a,v]_P \cap [b,w]_P$ the {\em inner of $C$}, and we call $C$ a {\em proper 4-crown in $P$} if $\J_P(C) = \emptyset$. Otherwise, we call it an {\em improper 4-crown in $P$}. 
\end{definition}
Due to $[a,v]_P \cap [b,w]_P = [a,w]_P \cap [b,v]_P$, this definition is independent of the choice of the disjoint edges of $C$. As indicated by the addendum ``in $P$'', a 4-crown is proper or improper only with respect to the poset it belongs to, not as an induced sub-poset as such. However, because the ambient poset is always clear in what follows, we omit explicit reference to it in most cases. In parts of the early literature \cite{Duffus_Rival_1979,Duffus_etal_1980_DRS,Rival_1976,Rival_1982,Rutkowski_1989}, improper 4-crowns are not regarded as crowns.

Examples for the different types of 4-crowns are shown in Figure \ref{fig_Intro}. In all three posets, the points $a,b,v,w$ form a 4-crown in $E(P)$. It is improper in the first two posets. In the third poset (it is Rutkowski's \cite{Rutkowski_1989} poset P9), the points $a,b,v,w$ form a proper 4-crown, but all edges of all crowns in $E(P)$ belong to improper 4-crowns.

From the different ways to construct new posets from old ones, we need the {\em direct sum} and the {\em ordinal sum} of two posets $P$ and $Q$ with disjoint carriers $X$ and $Y$. The direct sum $P + Q$ has the partial order relation ${\leq_P} \cup {\leq_Q}$, whereas the ordinal sum $P \oplus Q$ has the partial order relation ${\leq_P} \cup {\leq_Q} \cup (X \times Y)$. In particular, if $P$ and $Q$ are disjoint antichains, the edge-set of $P \oplus Q$ is simply $X \times Y$.

For a homomorphism $f$ between posets or graphs, we denote by $f \vert_Z$ the pre-restriction of $f$ to a subset $Z$ of its domain. An order homomorphism $r : P \rightarrow P$ is called a {\em retraction of the poset $P$} if $r$ is idempotent, and an induced sub-poset $R$ of $P$ is called a {\em retract of $P$} if a retraction $r : P \rightarrow P$ exists with $r[X] = R$. For the sake of simplicity, we identify $r$ with its post-restriction and write $r : P \rightarrow R$. A poset $P$ has the fixed point property iff every retract of $P$ has the fixed point property \cite[Theorem 4.8]{Schroeder_2016}. An improper 4-crown cannot be a retract of $P$ because its inner points cannot be mapped to it in an order preserving way.

Finally, we want to note that according to Duffus et al.\ \cite[p.\ 232]{Duffus_etal_1980_DPR}, the retracts of $P$ in $E(P)$ contain all information about retracts of height one in $P$: Given a retraction $r : P \rightarrow R$, there exists a retraction $r' : P \rightarrow Q$ with $Q \simeq R$ and $E(Q) \subseteq E(P)$. In particular, we have $Q \subseteq E(P)$ if $R$ has height one.

\subsection{Results about homomorphisms} \label{subsec_homs}

In this technical section, we compile required tools and results about homomorphisms. Here as in all following sections, $P = (X, \leq_P)$ and $Q = (Y, \leq_Q)$ are finite posets with at least two points and without isolated points. (For the existence of homomorphisms between posets and in particular of retractions, isolated points do not play a significant role.) The sets $L(P)$, $U(P)$ and $L(Q)$, $U(Q)$ are thus disjoint antichains in $P$ and $Q$.

\begin{lemma} \label{lemma_f_strict}
Let $f : P \rightarrow Q$ be a surjective homomorphism. There exists a homomorphism $g : P \rightarrow Q$ with $g[L(P)] = L(Q)$, $g[U(P)] = U(Q)$, and $g \vert_{M(P)} = f \vert_{M(P)}$. If $f$ is a retraction, also $g$ can be chosen as retraction.
\end{lemma}
\BP For every $y \in Y$, we select $\lambda(y) \in L(Q) \cap \darr_Q y$ and $\yps(y) \in U(Q) \cap \uarr_Q y $. We define a mapping $g : X \rightarrow Y$ by
\begin{align*} 
x & \mapsto
\begin{cases}
\lambda(f(x)), & \mytext{if } x \in L(P) \mytext{ and } f(x) \notin L(Q) , \\
\yps(f(x)), & \mytext{if } x \in U(P) \mytext{ and } f(x) \notin U(Q), \\
f(x), & \mytext{otherwise.}
\end{cases}
\end{align*}
Because $f$ is onto, we must have $L(Q) \subseteq f[L(P)]$. For $y \in L(Q)$, there exists thus an $x \in L(P)$ with $f(x) = y \in L(Q)$. Therefore, $g(x) = f(x) = y$, thus $y \in g[L(P)]$, and $g[L(P)] = L(Q)$ is shown because $g[L(P)] \subseteq L(Q)$ is clear. The equation $g[U(P)] = U(Q)$ is dual.

Let $x, y \in X$ with $x <_P y$. Then $x \notin U(P)$, thus $g(x) \leq_Q f(x)$, and $y \notin L(P)$, hence $f(y) \leq_Q g(y)$, and $g$ is a homomorphism.

Now assume that $f$ is a retraction. Then $Q$ is an induced sub-poset of $P$, and we have $L(P) \cap Q \subseteq L(Q)$. For $x \in L(P) \cap Q$, we thus have $f(x) = x \in L(Q)$, hence $g(x) = f(x)$. In the same way we see $g(x) = f(x)$ for all $x \in U(P) \cap Q$. For $x \in M(P) \cap Q$, $g(x) = f(x)$ is trivial, and $g \vert_Q = f \vert_Q = \id_Q$ is shown.

\EP 



The following proposition has already been proven in \cite[Proposition 1]{aCampo_2024}:

\begin{proposition} \label{prop_fortsetzung}
Let $Q$ be a poset of height one and let $f : E(P) \rightarrow Q$ be a homomorphism. We define for every $x \in X$
\begin{align*}
\alpha(x) & := f \left[ L(P) \cap \darr_P x \right], \quad
\beta(x)  := f \left[ U(P) \cap \uarr_P x \right].
\end{align*}
If, for all $x \in X$ with $\alpha(x) \cap \beta(x) = \emptyset$,
\begin{align} \label{alphabeta_cond}
\begin{split}
\# \alpha(x) \geq 2 \quad & \Rightarrow \quad \# \beta(x) = 1, \\
\# \beta(x) \geq 2 \quad & \Rightarrow \quad \# \alpha(x) = 1,
\end{split}
\end{align}
then there exists a homomorphism $g : P \rightarrow Q$ with $g \vert_{E(P)} = f$.
\end{proposition}

The next result is a variant of \cite[Lemma 1]{aCampo_2024}:

\begin{lemma} \label{lemma_fourCrown}
Let $Q$ be a poset of height one and let  $f : E(P) \rightarrow Q$ be a strict surjective homomorphism. $f$ can be extended to $P$ iff the following condition holds:
\begin{align} \label{bedingung_Part}
\begin{split}
& \mytext{A 4-element set } F \subseteq E(P) \mytext{ is not an improper 4-crown if it contains} \\
& \mytext{points from two different sets } f^{-1}(a), f^{-1}(b) \mytext{ with } a, b \in L(Q) \\
& \mytext{and from two different sets in } f^{-1}(v) f^{-1}(w) \mytext{ with } v, w \in U(Q).
\end{split} 
\end{align}
\end{lemma}

\BP  Only direction ``$\Leftarrow$'' has to be shown. For $x \in X$, we define the sets $\alpha(x)$ and $\beta(x)$ as in Proposition \ref{prop_fortsetzung}. Then $\alpha(x) \subseteq L(Q)$ and $\beta(x) \subseteq U(Q)$. Assume that $\alpha(x)$ contains two different points $a, b$ and $\beta(x)$ contains two different points $v, w$. Then $x \in M(P)$, and there exist points $a' \in \urbild{f}(a) \cap \darr_P x$, $b' \in \urbild{f}(b) \cap \darr_P x$, $v' \in \urbild{f}(v) \cap \uarr_P x$, and $w' \in \urbild{f}(w) \cap \uarr_P x$. But then $F := \setx{ a', b', v', w'} \subseteq E(P)$ is an improper 4-crown of $P$ containing points from $\urbild{f}(a), \urbild{f}(b)$ and $\urbild{f}(v), \urbild{f}(w)$ in contradiction to \eqref{bedingung_Part}. At least one of the sets $\alpha(x)$ and $\beta(x)$ is thus a singleton, and Proposition \ref{prop_fortsetzung} yields an extension $g : P \rightarrow Q$ of $f$.

\EP

\section{Posets of height one as retracts} \label{sec_FourCrowns}

\subsection{The multigraphs $\fC(C)$ and $\fFP$} \label{subsec_Cgraph}

\begin{figure}
\begin{center}
\includegraphics[trim = 70 630 295 74, clip]{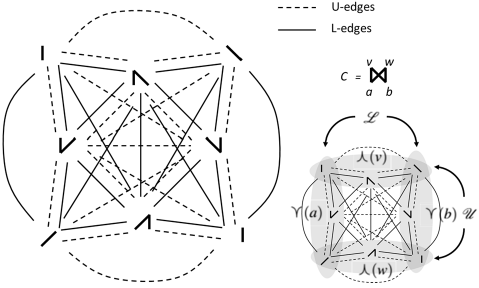}
\caption{\label{fig_CGraph} The multigraph $\fC(C)$ with loops omitted in the case of a 4-crown $C$. Explanation in text.}
\end{center}
\end{figure}

In this as in the following sections, the poset $C$ is a connected poset of height one with carrier $Z$. We thus have $Z = L(C) \cup U(C)$ and ${<_C} \subseteq L(C) \times U(C)$.

In Lemma \ref{lemma_fourCrown}, we have identified  the points of improper 4-crowns in $E(P)$ as the objects crucial for the existence of an extension of an homomorphism from $E(P)$ onto $C$. However, it is beneficial to widen the scope from improper 4-crowns to objects comprising the points of all improper 4-crowns having a common inner point:

\begin{definition} \label{def_crownBundle}
We define $\C_{\imp}(P)$ as the set of improper 4-crowns contained in $E(P)$. Furthermore, we define
\begin{align*}
\MP := & \cup_{F \in \C_{\imp}(P)} \J_P(F) \mytext{ (the set of all inner points of improper 4-crowns),} \\
\Xi(m) := & \left( L(P) \cap \darr m \right) \; \cup \; \left( U(P) \cap \uarr m \right) \quad \mytext{for all } m \in \MP, \\
& \mytext{(the union of all improper 4-crowns with inner point m),} \\
\F_0(P) := & \mysetdescr{ \Xi(m) }{ m \in \MP }, \\
\FP := & \mytext{the set of maximal sets in } \F_0(P) \mytext{ with respect to set inclusion.}
\end{align*}
The elements of $\FP$ are called {\em 4-crown bundles}.
\end{definition}

For an edge $x <_{E(P)} y$ in $E(P)$, the existence of a 4-crown bundle containing $x$ and $y$ is equivalent to the existence of an improper 4-crown containing these points.

The image of a set $\Xi(m)$ under a homomorphism from $P$ to $C$ cannot contain several points of $L(C)$ {\em and} several points of $U(C)$. The following sets contain thus all possible non-singleton images of 4-crown bundles under order-preserving mappings to $L(C) \oplus U(C)$: for every $a \in L(C)$ and every $v \in U(C)$
\begin{align*}
\Mvee(a) & \; := \; \mysetdescr{ \setx{a} \cup N }{ \emptyset \not= N \subseteq U(C) }, \quad \L := \bigcup_{a \in L(C)} \Mvee(a), \\
\Mwedge(v) & \; := \; \mysetdescr{ \setx{v} \cup N }{ \emptyset \not= N \subseteq L(C) }, \quad \; \U := \bigcup_{v \in U(C)} \Mwedge(v).
\end{align*}

The next three definitions describe the graphs and the type of homomorphisms we will work with. In the first one, we select the elements of $\L \cup \U$ belonging to $C$ and make them a multigraph:

\begin{definition} \label{def_Cgraph}
Let
\begin{align*}
\Co(C) := \quad & \mysetdescr{S \in \L }{ \exists \; a \in L(C) \mytext{ : } S \subseteq \uarr_C \; a } \\
\cup & \mysetdescr{S \in \U }{ \exists \; v \in U(C) \mytext{ : } S \subseteq \darr_C \; v }.
\end{align*}
We define $\fC(C)$ as the undirected multigraph with vertex set $\Co(C)$ and an edge set containing two types of edges, {\em L-edges} and {\em U-edges}, which are assigned as follows: For $S, T \in \Co(C)$, the edge set of $\fC(C)$ contains
\begin{itemize}
\item an  L-edge between $S$ and $T$ iff $L(C) \cap S \cap T \not= \emptyset$;
\item a U-edge between $S$ and $T$ iff $U(C) \cap S \cap T \not= \emptyset$.
\end{itemize}
We call the two-element sets in $\Co(C)$ the {\em tips} of $\fC(C)$. 
\end{definition}
For a 4-crown $C$, the multigraph $\fC(C)$ with loops omitted is shown in Figure \ref{fig_CGraph}.

A tip contains the points belonging to a single edge of $C$. It is thus contained in exactly one of the sets $\Mvee(a)$ {\em and} in exactly one of the sets $\Mwedge(v)$, whereas each $S \in \Co(C)$ not being a tip is contained {\em either} in exactly one of the sets $\Mvee(a)$ {\em or} in exactly one of the sets $\Mwedge(v)$. We note that no L-edge exists between vertices belonging to $\Mvee(a)$ and $\Mvee(b)$ with $a \not= b$. Correspondingly, vertices in $\Mwedge(v)$ and $\Mwedge(w)$ with $v \not= w$ are not connected by an U-edge. 

According to Lemma \ref{lemma_fourCrown}, the existence of an homomorphism $g$ from $P$ onto $C$ with $g \vert_{E(P)}$ being strict depends entirely on the possibility to map the points of 4-crown bundles in $E(P)$ to $C$ in an appropriate way. Every 4-crown bundle in $E(P)$ is mapped to an element of $\Co(C)$, and for $F$ and $G$ being 4-crown bundles in $E(P)$ with $F \cap G \not= \emptyset$, also their images $g[F], g[G] \in \Co(C)$ have to share points. The complexity of the network of 4-crown bundles in $P$ is going to be reflected by their images $g[F]$ in $\fC(C)$.

\begin{definition} \label{def_fFP}
The graph $\fFP$ is the undirected multigraph with vertex set $\FP$ and an edge set containing two types of edges called {\em L-edges} and {\em U-edges} again, which are assigned as follows: For $F, G \in \FP$, the edge set of $\fFP$ contains
\begin{itemize}
\item an L-edge between $F$ and $G$ iff $L(P) \cap F \cap G \not= \emptyset$;
\item a U-edge between $F$ and $G$ iff $U(P) \cap F \cap G \not= \emptyset$.
\end{itemize}
\end{definition}


\begin{definition} \label{def_LUhom}
We call a homomorphism $\phi : \fFP \rightarrow \fC(C)$ {\em separating} if there exists a one-to-one mappings $i : Z \rightarrow E(P)$ with $i[L(C)] \subseteq L(P)$ and $i[U(C)] \subseteq U(P)$ providing the following property: for all $a \in L(C)$, $v \in U(C)$ and all $F \in \FP$
\begin{align} \label{phi_cond}
\begin{split}
i(a) \in F \mytext{ and } \phi(F) \in \L & \quad \Rightarrow \quad \phi(F) \in \Mvee(a), \\
i(v) \in F \mytext{ and } \phi(F) \in \U & \quad \Rightarrow \quad \phi(F) \in \Mwedge(v).
\end{split}
\end{align}
We say that $i$ {\em belongs to $\phi$.} Furthermore, we call a separating homomorphism $\phi$ {\em $Z$-separating} in the case of $Z \subseteq E(P)$ and $i[Z] = Z$.
\end{definition}

If there exists a separating homomorphism $\phi$, the set $i[Z] \subseteq E(P)$ cannot contain an improper 4-crown $\setx{i(a),i(b)} \oplus \setx{i(v),i(w)}$. Otherwise, there exists an $F \in \FP$ with $\setx{ i(a), i(b), i(v), i(w) } \subseteq F$. In the case of $\phi(F) \in \L$, \eqref{phi_cond} yields $\phi(F) \in \Mvee(a) \cap \Mvee(b)$, but this intersection is empty. The case $\phi(F) \in \U$ is analogous. Moreover, for a 4-crown bundle $F \in \FP$ containing more than one point of $i[L(C)]$ (more than one point of $i[U(C)]$), the image $\phi(F)$ must belong to $\U$ (to $\L$).

\subsection{Separating homomorphisms and retracts} \label{subsec_sepHomsAndRetracts}

The next two propositions illuminate the interdependence between surjective homomorphisms from $P$ onto $C$ and separating homomorphisms from $\fFP$ to $\fC(C)$. One of the dependencies is straight forward:

\begin{lemma} \label{lemma_surjAufC}
The existence of a homomorphism from $P$ onto $C$ implies the existence of a separating homomorphism from $\fFP$ to $\fC(C)$ with $i^{-1} : E(P) \vert_{i[Z]} \rightarrow C$ being a homomorphism. Moreover, $C \subseteq E(P)$ being a retract of $P$ implies the existence of a $Z$-separating homomorphism from $\fFP$ to $\fC(C)$ with $i$ being an automorphism of $C$.
\end{lemma}
\BP Let $f : P \rightarrow C$ be a surjective homomorphism. According to Lemma \ref{lemma_f_strict}, there exists a surjective homomorphism $g : P \rightarrow C$ with $g[L(P)] = L(C)$ and $g[U(P)] = U(C)$. The mapping
\begin{align*} 
\begin{split}
\phi : \FP & \rightarrow \Co(C) \\
F & \mapsto g[ F ],
\end{split}
\end{align*}
is thus a homomorphism from $\fFP$ to $\fC(C)$. For every $a \in L(C)$, we select a point $i(a) \in g^{-1}(a)$. For $F \in \FP$ with $i(a) \in F$, the membership $\phi(F) \in \L$ means that $a = g(i(a)) \in \phi(F)$ is the only point of $L(C)$ contained in $g[F]$, thus $\phi(F) \in \Mvee(a)$. The argumentation for $v \in U(C)$ is analogous, and $\phi$ is separating. $i(x) <_P i(y)$ implies $x = g(i(x)) <_C g(i(y)) = y$ for all $x, y \in Z$, and $i^{-1}$ is a homomorphism. 

Let $C \subseteq E(P)$. For $f : P \rightarrow C$ being a retraction, the addendum in Lemma \ref{lemma_f_strict} states that also our homomorphism $g$ can be chosen as retraction. Now select $i := \id_Z$.

\EP

For the other dependency, we make the assumption that $C$ shares some features with an ordinal sum of antichains:


\begin{lemma} \label{lemma_heightOne}
Assume that there exist two non-empty sets $L_0 \subseteq L(C)$ and $U_0 \subseteq U(C)$ with $(L_0 \times U(C)) \cup (L(C) \times U_0) \subseteq {<_C}$. There exists a surjective homomorphism from $P$ onto $C$ if a separating homomorphism $\phi : \fFP \rightarrow \fC(C)$ exists with
\begin{align} \label{Ann_L0}
& \mysetdescr{ a \in L(C) }{ \phi[ \FP ] \cap \Mvee(a) \not= \emptyset } \subseteq L_0, \\ \label{Ann_U0}
& \mysetdescr{ v \in U(C) }{ \phi[ \FP ] \cap \Mwedge(v) \not= \emptyset } \subseteq U_0, \\ \nonumber
& \; i^{-1} : E(P) \vert_{i[C]} \rightarrow C \mytext{ is a homomorphism,}
\end{align}
where $i : Z \rightarrow E(P)$ is the one-to-one mapping belonging to $\phi$. Moreover, in the case of $C \subseteq E(P)$, $C$ is a retract of $P$ if a $Z$-separating homomorphism $\phi$ fulfills \eqref{Ann_L0} and \eqref{Ann_U0} with $i$ being an automorphism of $C$.

\end{lemma}
\BP Let
\begin{align*}
\L_0 & := \mysetdescr{ a \in L(C) }{ \phi[ \FP ] \cap \Mvee(a) \not= \emptyset }, \\
\U_0 & := \mysetdescr{ v \in U(C) }{ \phi[ \FP ] \cap \Mwedge(v) \not= \emptyset }.
\end{align*}
We define for every $a \in L(C)$
\begin{align*}
V_0(a) & \; := \; \mysetdescr{ x \in L(P)}{ \exists F \in \FP \mytext{:} x \in F \mytext{ and } \phi(F) \in \Mvee(a) }
\end{align*}
Let $a, b \in L(C)$ with $a \not= b$ and let $x \in V_0(a)$, $y \in V_0(b)$. There exist $F, G \in \FP$ with $x \in L(P) \cap F$, $\phi(F) \in \Mvee(a)$, and $y \in L(P) \cap G$, $\phi(G) \in \Mvee(b)$. Due to $a \not= b$, there exists no L-edge between $\phi(F)$ and $\phi(G)$ in $\fC(C)$, and because $\phi$ is a homomorphism, we must have $x \not= y$. The sets $V_0(a)$, $a \in L(C)$, are thus pairwise disjoint (but can be empty).

Let $a \in L(C)$ with $i(a) \notin V_0(a)$. According to \eqref{phi_cond}, there exists no $F \in \FP$ with $i(a) \in F$ and $\phi(F) \in \L$. $\phi(F) \in \L$ for an $F \in \FP$ implies thus $i(a) \notin F$, hence $i(a) \notin \cup_{b \in L(C)} V_0(b)$. Setting
\begin{align*}
\forall \; a \in L(C) \mytext{ :} \quad V_1(a) & :=
\begin{cases}
V_0(a), & \mytext{if } i(a) \in V_0(a), \\
V_0(a) \cup \setx{i(a)} & \mytext{otherwise,}
\end{cases}
\end{align*}
yields thus pairwise disjoint non-empty subsets of $L(P)$ with $i(a) \in V_1(a)$ for all $a \in L(C)$ and $V_1(a) = \setx{i(a)}$ for all $a \in L(C) \setminus \L_0$.

Selecting an arbitrary $\ell \in L_0$, we define
\begin{align*}
\L_1 & :=
\begin{cases}
\L_0, & \mytext{if } \L_0 \not= \emptyset, \\
\setx{ \ell }, & \mytext{otherwise,}
\end{cases}
\end{align*}
and we distribute the elements of $L(P) \setminus \cup_{a \in L(C)} V_1(a)$ arbitrarily on the sets $V_1(a)$ with $a \in \L_1$. The resulting sets $V_2(a) \supseteq V_1(a)$ form a partition of $L(P)$ with the following properties
\begin{align*}
\forall \; a \in L(C) \mytext{ :} \quad i(a) & \in V_2(a), \\
\forall \; a \in L(C) \setminus \L_1 \mytext{ :} \quad V_2(a) & = \setx{i(a)}.
\end{align*}

In the same way, we treat the points of $U(P)$. 
With an arbitrary $u \in U_0$, we define
\begin{align*}
\forall \; v \in U(C) \mytext{ :} \quad \Lambda_0(v) & := \mysetdescr{ x \in U(P)}{ \exists F \in \FP \mytext{:} x \in F \mytext{ and } \phi(F) \in \Mwedge(v) }, \\
\forall \; v \in U(C) \mytext{ :} \quad \Lambda_1(v) & :=
\begin{cases}
\Lambda_0(v), & \mytext{if } i(v) \in \Lambda_0(v), \\
\Lambda_0(v) \cup \setx{i(v)} & \mytext{otherwise,}
\end{cases} \\
\U_1 & :=
\begin{cases}
\U_0, & \mytext{if } \U_0 \not= \emptyset, \\
\setx{ u }, & \mytext{otherwise,}
\end{cases}
\end{align*}
and we distribute the elements of $U(P) \setminus \cup_{v \in U(C)} \Lambda_1(v)$ arbitrarily on the sets $\Lambda_1(v)$ with $v \in \U_1$. The resulting sets $\Lambda_2(v) \supseteq \Lambda_1(v)$ form a partition of $U(P)$ fulfilling
\begin{align*}
\forall \; v \in U(C) \mytext{ :} \quad i(v) & \in \Lambda_2(v), \\
\forall \; v \in U(C) \setminus \U_1 \mytext{ :} \quad \Lambda_2(v) & = \setx{i(v)}.
\end{align*}

Now we define a mapping $f$ from $E(P)$ onto $Z$ by setting 
\begin{align*}
f^{-1}(a) & :=  V_2(a) \quad \mytext{for all } a \in L(C), \\
f^{-1}(v) & := \Lambda_2(v) \quad \mytext{for all } v \in U(C).
\end{align*}

Firstly, we show that $f$ is a homomorphism. Let $x \in L(P)$ and $y \in U(P)$ with $f(x) \parallel_C f(y)$. Because of 
$( \L_1 \times U(C) ) \cup ( L(C) \times \U_1 ) \subseteq {<_C}$, we have $f(x) \notin \L_1$, $f(y) \notin \U_1$, thus $x \in f^{-1}(f(x)) = V_2( f(x) ) = \setx{ i(f(x)) }$ and similarly $y \in \setx{ i(f(y)) }$. Because $i^{-1} : E(P) \vert_{i[C]} \rightarrow C$ is a homomorphism, we must have $x \parallel_{E(P)} y$, as desired.

We finish the proof by showing that $f$ fulfills Condition \eqref{bedingung_Part} in Lemma \ref{lemma_fourCrown}. In the case of $\L_0 = \emptyset = \U_0$, the set $\F(P)$ is empty. There exists thus no improper 4-crown in $E(P)$, and the condition is trivially fulfilled. Otherwise, let $G$ be an improper 4-crown and let $F \in \FP$ be a 4-crown bundle with $G \subseteq F$. In the case of $\phi(F) \in \Mvee(a)$ for an $a \in L(C)$, we have $L(P) \cap F \subseteq V_0(a)$, thus $G \cap f^{-1}(b) \subseteq F \cap f^{-1}(b) = F \cap V_2(b) = \emptyset$ for all $b \in L(C) \setminus \setx{a}$. In the same way, we see that for $\phi(F) \in \Mwedge(v)$ for a $v \in U(C)$, we must have $G \cap f^{-1}(w) = \emptyset$ for all $w \in U(C) \setminus \setx{v}$. Condition \eqref{bedingung_Part} is thus fulfilled in all cases.

For the addendum, observe $z = f(i(z))$ for all $z \in Z$.




\EP

In the light of these results, the reader may ask why we did not postulate in Definition \ref{def_LUhom} that the inverse of an one-to-one mapping belonging to a separating homomorphism has to be order-preserving. The reason is Proposition \ref{prop_clique} in Section \ref{subsec_imagePhi} where we must allow also other one-to-one mappings.

Our main result about ordinal sums of height one is:
\begin{theorem} \label{theo_surjAufC}
Let $C$ be an ordinal sum of height one. The existence of a surjective homomorphism from $P$ onto $C$ is equivalent to the existence of a separating homomorphism from $\fFP$ to $\fC(C)$. Moreover, if $C$ is a subposet of $E(P)$, then $C$ is a retract of $P$ iff a $Z$-separating homomorphism from $\fFP$ to $\fC$ exists.
\end{theorem}
\BP For both propositions, Lemma \ref{lemma_surjAufC} delivers direction ``$\Rightarrow$''. Assume that a separating homomorphism $\phi$ from $P$ to $\fC$ exists and let $i$ be the one-to-one mapping belonging to $\phi$. Because of $C = L(C) \oplus U(P)$, the inverse mapping $i^{-1}$ is a homomorphism, and with $L_0 := L(C)$, $U_0 := U(C)$, all assumptions in Lemma \ref{lemma_heightOne} are fulfilled. In the case that $C$ is a subposet of $E(P)$ and $\phi$ is $Z$-separating, $i$ is an automorphism of $C$, and the proposition follows with the addendum in Lemma \ref{lemma_heightOne}.

\EP

If $P$ is a poset of height one, then $\FP = \emptyset$, and for every $C \subseteq E(P)$ being an ordinal sum of antichains, the trivial mapping $(\emptyset, \emptyset, \Co(C))$ from $\FP$ to $\Co(C)$ is a $Z$-separating homomorphism from $\fFP = (\emptyset,\emptyset)$ to $\fC(C)$. Theorem \ref{theo_surjAufC} states that $C$ is a retract of $P$ and confirms thus the well-known result of Rival \cite{Rival_1976,Rival_1982} that a 4-crown in a poset of height one is always a retract.

The mapping $(\emptyset, \emptyset, \Co(C))$ is an extreme example for the fact that a separating homomorphism $\phi : \fFP \rightarrow \fC(C)$ must not tell much about a corresponding (strict) surjective homomorphism $f : E(P) \rightarrow C$. Tracing back the construction process in the proof of Lemma \ref{lemma_fourCrown}, we realize that $f$ has to send the points contained in the sets $V_1(v)$ and $\Lambda_1(a)$ to different points in $L(C)$ and $U(C)$, respectively. With respect to the images of the remaining points in $L(P)$ and $U(P)$, each distribution onto the final sets $V_2(a)$ and $\Lambda_2(v)$ with $a \in \L_1$, $v \in \U_1$ is possible.

Let $C = L(C) \oplus U(C)$, $Y := \left( \cup \; \FP \right) \cup \MP \cup C$, and $Q := P \vert_Y$. Because of $\fFP = \fFQ$ and Theorem \ref{theo_surjAufC}, the poset $C$ is a retract of $P$ iff it is a retract of $Q$. The elements of $X \setminus Y$ are thus irrelevant for $C$ being a retract of $P$ or not. And indeed, for $x \in M(P) \setminus Y$, one of the sets $L(P) \cap \darr x $ or $U(P) \cap \uarr x$ is a singleton, and according to \cite[Ex.\ 4-5]{Schroeder_2016}, the point $x$ can be removed by a series of I-retractions. By doing so for all points in $M(P) \setminus Y$, we arrive at the poset $P' := P \vert_{E(P) \cup Y}$. The poset $P' \vert_Y$ contains $C$, the improper 4-crowns of $P$ and their inner, and the poset $P' \vert_{E(P) \setminus Y}$ has height one or is an antichain. Intuitively, it is plausible that a homomorphism mapping $P' \vert_Y$ onto $C$ can be extended to the total poset $P'$.

\begin{figure}
\begin{center}
\includegraphics[trim = 75 635 165 70, clip]{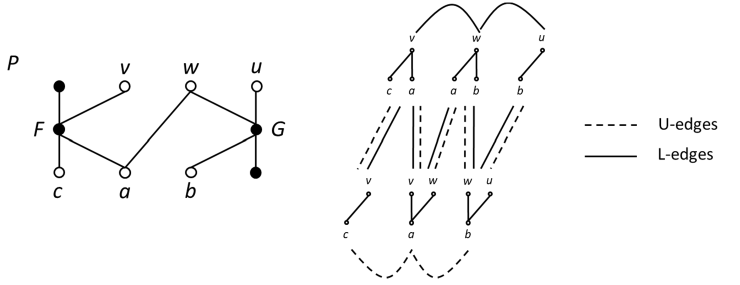}
\caption{\label{fig_Cmax} Left: A poset $P$ containing two 4-crown bundles $F$ and $G$. Right: The graph $\fCdr(C)$ for the subposet $C \subset E(P)$ marked by hollow dots in $P$. The lower row of vertices in $\fCdr(C)$ is from left to right $\vee_C(c)$, $\vee_C(a)$, and $\vee_C(b)$, the upper row is $\wedge_C(v)$, $\wedge_C(w)$, and $\wedge_C(u)$.}
\end{center}
\end{figure}

Let $m := \# L(C)$, $n := \# U(C)$. Containing up to $ n ( 2^m - 1) + m (2^n -1) - mn$ vertices, the graph $\fC(C)$ quickly becomes unwieldy even for small values of $m$ and $n$. However, we can replace it by a considerably simpler target-graph:

\begin{corollary} \label{coro_Cdrei}
Let for every $a \in L(C)$, $v \in U(C)$
$$
\vee_C(a) := \uarr_C \; a, \quad \wedge_C(v) := \darr_C \; v,
$$
and define
\begin{align*}
\Codr(C) & \; := \; \mysetdescr{ \vee_C(a) }{ a \in L(C) } \; \cup \; \mysetdescr{ \wedge_C(v) }{ v \in U(C) }, \\
\fCdr(C) & \; := \; \mytext{the sub-graph of } \fC \mytext{ induced by} \Codr(C).
\end{align*}
There exists a separating ($Z$-separating) homomorphism from $\fFP$ to $\fC(C)$ iff a separating ($Z$-separating) homomorphism from $\fFP$ to $\fCdr(C)$ exists.
\end{corollary}
\BP Only direction ``$\Rightarrow$'' has to be shown. For the proof, we define the mapping $\Theta : \Co(C) \rightarrow \Codr(C)$ as follows:
\begin{align*}
\Theta(S) & :=
\begin{cases}
\vee_C(a), & \mytext{if } S \in \Mvee(a), \\
\wedge_C(v), & \mytext{if } S \in \Mwedge(v) \mytext{ is not a tip.}
\end{cases}
\end{align*}
Let $\phi$ be a separating ($Z$-separating) homomorphism from $\fFP$ to $\fC(C)$. Because of $S \subseteq \Theta(S)$ for all $S \in \Co(C)$, the mapping $\Theta \circ \phi$ is a homomorphism from $\fFP$ to $\fCdr(C)$.

Let $F \in \FP$. If $\Theta(\phi(F))$ is not a tip, then $\Theta(\phi(F)) \in \N$ implies $\phi(F) \in \N$ for both $\N \in \setx{ \L, \U }$, and \eqref{phi_cond} is fulfilled. And if $\Theta(\phi(F))$ is a tip, then $\phi(F) \subseteq \Theta(\phi(F))$ yields $\phi(F) = \Theta(\phi(F))$, and \eqref{phi_cond} holds in this case, too.

\EP

Figure \ref{fig_Cmax} presents an example. On the left, a poset $P$ is shown containing two 4-crown bundles $F$ and $G$. For the subposet $C \subset E(P)$ marked by hollow dots, the graph $\fCdr(C)$ is shown on the right.

\subsection{Cliques as images of homomorphisms} \label{subsec_imagePhi}

We call a non-empty subset $\K$ of $\FP$ or of $\Co(C)$ a {\em clique} if the vertices contained in $\K$ are all pairwise connected by an L-edge {\em and} a U-edge. Because each vertex in $\fC(C)$ has an L-loop and a U-loop, every mapping from $\fFP$ to a clique in $\fC(C)$ is a homomorphism. In checking if a given poset $P$ has $C \subseteq E(P)$ as retract, it may thus be tempting to look how $\FP$ can be mapped to a clique in $\fC(C)$. We see in this section that this is possible in a $Z$-separating way iff the improper 4-crowns in $P$ contain only a sparse subset of the edges of $C$.

\begin{figure}
\begin{center}
\includegraphics[trim = 70 705 365 74, clip]{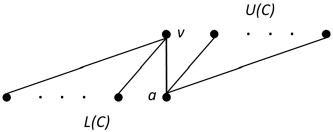}
\caption{\label{fig_Nav} The poset $N_{a,v}$ contains $m + n - 1$ edges of the up to $m n$ edges of $C$.}
\end{center}
\end{figure}

Let thus $C \subseteq E(P)$. We need two posets with carrier $Z$:
\begin{itemize}
\item The poset $C_{\FP}$ has carrier $Z$ and edge set
$$
<_{C_{\FP}} := \mysetdescr{ (a,v) \in {<_C} }{ \exists \; F \in \FP \mytext{ : } a,v \in F }.
$$
\item For $a \in L(C)$, $v \in U(C)$, the poset $N_{a,v}$ has carrier $Z$ and edge set
$$
<_{N_{a,v}} := ( L(C) \times \setx{v} ) \cup ( \setx{a} \times U(C) ).
$$
\end{itemize}

$N_{a,v}$ is illustrated in Figure \ref{fig_Nav}. Our result is:
\begin{proposition} \label{prop_clique}
Let $C \subseteq E(P)$. There exists a $Z$-separating homomorphism $\phi$ from $\fFP$ to $\fC(C)$ for which the image set $\phi[ \FP ]$ is a clique iff there exist points $a \in L(C)$, $v \in U(C)$ with $C_{\FP}$ being a subposet of $N_{a,v}$. 
\end{proposition}
\BP ``$\Rightarrow$'': Let $\phi: \fFP \rightarrow \fC(C)$ be a $Z$-separating homomorphism with $\phi[ \FP ]$ being a clique in $\fC(C)$. With the mapping $\Theta : \Co(C) \rightarrow \Codr(C)$ defined in the proof of Corollary \ref{coro_Cdrei}, the mapping $\psi := \Theta \circ \phi$ is a $Z$-separating homomorphism. The set $\K := \psi[ \FP ]$ is a clique in $\fC_{\max}(C)$ and has therefore cardinality zero, one, or two.

$\K = \emptyset$ yields $\FP = \emptyset$, thus $C_{\FP} = L(C) + U(C)$.

$\K = \setx{S}$ with $S \in \Codr(C)$: Assume $S \in \L$. Due to \eqref{phi_cond}, the set $L(C) \cap (\cup \FP) = L(C) \cap (\cup \psi^{-1}(S))$ contains at most a single point. If it is empty, then $C_{\FP}$ is the antichain $L(C) + U(C)$, and in the case of $L(C) \cap ( \cup \FP) = \setx{a}$, the poset $C_{\FP}$ is a subposet of $( L(C) \setminus \setx{a} ) + ( \setx{a} \oplus U(C) )$. The case $S \in \U$ is analogous.

$\K = \setx{S, T}$ with $S, T \in \Codr(C)$, $S \not= T$: Let $i$ be the bijection belonging to $\psi$. One of the sets $S, T$ must belong to $\L$ and the other one to $\U$. Let $S = \vee_C(a)$ and $T = \wedge_C(v)$ with $a \in L(C)$, $v \in U(C)$. Condition \eqref{phi_cond} enforces
\begin{align*}
\begin{split}
L(C) \cap ( \cup \psi^{-1}(S) ) & \; \subseteq \; \setx{ i(a)}, \\
U(C) \cap ( \cup \psi^{-1}(T) ) & \; \subseteq \; \setx{ i(v)}.
\end{split}
\end{align*}
Let $(x,y) \in {<_{C_{\FP}}}$ and let $F \in \FP$ be a 4-crown bundle with $x, y \in F$. Then $x,y \in \psi^{-1}( \psi(F))$, thus $x = i(a)$ or $y = i(v)$. We conclude that  $C_{\FP}$ is a subposet of $N_{i(a),i(v)}$.



``$\Leftarrow$'': Let $a \in L(C)$, $v \in U(C)$ with $C_{\FP}$ being a subposet of $N_{a,v}$. We distinguish four cases. In the first three ones, the bijection $i$ belonging to the constructed $Z$-separating homomorphism is simply $\id_Z$.

1. $L(C) = \setx{a}$ and $U(C) = \setx{v}$: We define $\phi(F) := \setx{a,v}$ for all $F \in \FP$.

2. $L(C) = \setx{a}$ and $U(C) \supset \setx{v}$: We have $\vee_C(a) = Z$ and a $Z$-separating homomorphism is given by $\phi(F) := \vee_C(a) \in \L \setminus \U$ for all $F \in \FP$.

3. $L(C) \supset \setx{a}$ and $U(C) = \setx{v}$: as the previous case.

4. $L(C) \supset \setx{a}$ and $U(C) \supset \setx{v}$: There exists an edge $(b,w) \in {<_C}$ with neither $\vee_C(b)$ nor $\wedge_C(w)$ being a tip. Because $C_{\FP}$ is a subposet of $N_{a,v}$, we have $C \cap F \subseteq L(C) \cup \setx{v}$ or $C \cap F \subseteq U(C) \cup \setx{a}$ for all $F \in \FP$. We define for all $F \in \FP$, $z \in Z$
\begin{align*}
\phi(F) & :=
\begin{cases}
\wedge_C(w), & \mytext{if } C \cap F \subseteq L(C) \cup \setx{v}, \\
\vee_C(b)& \mytext{otherwise,}
\end{cases} \\
i(z) & :=
\begin{cases}
b, & \mytext{if } z = a, \\
a, & \mytext{if } z = b, \\
w, & \mytext{if } z = v, \\
v, & \mytext{if } z = w, \\
z, & \mytext{otherwise.}
\end{cases}
\end{align*}
Let $F \in \FP$. Trivially, $a \in F$ and $\phi(F) \in \L$ yields $\phi(F) = \vee_C(b)$, and for $x \in F \cap L(C)$ with $x \not= a$, we have $C \cap F \not\subseteq  U(C) \cup \setx{a}$, thus $\phi(F) = \wedge_C(w) \notin \L$. The cases $v \in F$ and $y \in F \cap ( U(C) \setminus \setx{v} )$ are analogous, and $\phi$ is $Z$-separating with $\phi[\FP]$ being a clique.

\EP

\begin{figure}
\begin{center}
\includegraphics[trim = 75 710 210 74, clip]{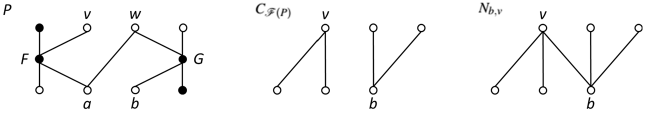}
\caption{\label{fig_Nav_Anm} A poset $P$ containing two 4-crown bundles $F$ and $G$. The subposet $C$ (hollow dots) is a retract and $C_{\FP}$ is a subposet of $N_{b,v}$.}
\end{center}
\end{figure}

The reader will have noticed that in the proofs of the third part of direction ``$\Rightarrow$'' and the fourth part of direction ``$\Leftarrow$'', neither the one-to-one mapping $i$ nor its inverse must be order-preserving.

Figure \ref{fig_Nav_Anm} illustrates the proposition. On the left, a poset $P$ containing two 4-crown bundles $F$ and $G$ is shown. (It is the poset already used in Figure \ref{fig_Cmax}.) The subposet $C$ marked by hollow dots is a retract of $P$. Because $F$ contains two points of $L(C)$ and $G$ two points of $U(C)$, every separating homomorphism must map $F$ into $\U$ and $G$ into $\L$, and setting $\phi(F) := \wedge_C(w) = \setx{a,b,w}$ and $\phi(G) := \vee_C(a) = \setx{a,v,w}$ yields a $Z$-separating homomorphismus with $\phi[\FP]$ being a clique. And indeed, $C_{\FP}$ is a subposet of $N_{b,v}$. We have $b \not<_C v$, and there exist no other points $x \in L(C)$, $y \in U(C)$ with $C_{\FP}$ being a subposet of $N_{x,y}$. 

For a 4-crown $C$, the posets $N_{a,v}$ with $a \in L(C)$, $v \in U(C)$, are exactly the four subposets of $C$ with $N$-shaped diagram. We conclude with Proposition \ref{prop_clique} and Theorem \ref{theo_surjAufC}:
\begin{corollary} \label{coro_clique}
Let $C \subseteq E(P)$ be a 4-crown. If $\fFP$ is a complete graph, then $C$ is a retract of $P$ iff there exists an edge of $C$ not contained in any 4-crown bundle.
\end{corollary}

\begin{figure}
\begin{center}
\includegraphics[trim = 70 719 290 70, clip]{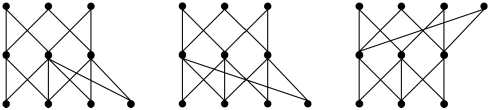}
\caption{\label{fig_P345} Rutkowski's posets P3, P4, and P5.}
\end{center}
\end{figure}

By means of this corollary we see that the posets P3-P5 of Rutkowski \cite{Rutkowski_1989} shown in Figure \ref{fig_P345} do not have a 4-crown as retract. The graph $\fFP$ of each of them is a complete graph with three vertices, but all edges of 4-crowns contained in their extremal points belong to improper 4-crowns.

A consequence of the corollary is that a poset $P$ with at most three maximal points and at most three minimal points has a 4-crown $C \subseteq E(P)$ as retract iff one of the edges of $C$ is not contained in an improper 4-crown. With this criterion, we see immediately that also the posets P1, P2, and P6-P10 of Rutkowski \cite{Rutkowski_1989} do not have a 4-crown as retract. (P9 is shown in Figure \ref{fig_Intro} on the right.)

\end{document}